\documentclass[a4paper,11pt]{article}
\usepackage{amssymb}

\newcommand{\bc}{\begin{center}}
\newcommand{\ec}{\end{center}}
\newcommand{\be}{\begin{equation}}
\newcommand{\ee}{\end{equation}}
\newcommand{\bea}{\begin{eqnarray}}
\newcommand{\eea}{\end{eqnarray}}
\newcommand{\ba}{\begin{array}}
\newcommand{\ea}{\end{array}}

\newcommand{\dsf}{\displaystyle\frac}

\def\l{\lambda}
\def\g{\gamma}
\def\r{\rho}
\def\d{\delta}

\def\m{\mu}

\def\i{\varepsilon}
\def\Q{\mathbb{Q}}
\def\R{\mathbb{R}}
\def\Z{\mathbb{Z}}
\def\N{\mathbb{N}}
\def\C{\mathbb{C}}
\begin{document}
\bc
{\Large {\bf ON RATIONAL $P$-ADIC DYANAMICAL SYSTEMS}}\\[6mm]

{\large FARRUKH MUKHAMEDOV \ \  UTKIR ROZIKOV}\\[4mm]
\ec 
\begin{abstract}
In the paper we investigate the behavior of trajectory of rational
$p$-adic dynamical system in complex $p$-adic filed $\C_p$. It is
studied Siegel disks and attractors of such dynamical systems. We
show that Siegel disks may either coincide or disjoin for
different fixed points of the dynamical system. Besides, we find
the basin of the attractor of the system. It is proved that such
kind of dynamical system is not ergodic on a unit sphere with
respect to the Haar measure.
\\[2mm]
{\it Key words and phrases:} Rational dynamics, attractors, Siegel
disk, complex $p$-adic field.\\
 {2000 {\it Mathematical Subject Classification:}} 46S10, 12J12, 11S99, 30D05, 54H20.
\end{abstract}
\normalsize

\section{Introduction}

The $p$-adic numbers were first introduced by the Germen
mathematician K.Hensel. For about a century after the discovery of
$p$-adic numbers, they were mainly considered objects of pure
mathematics. Beginning with 1980's various models described in the
language of $p$-adic analysis have been actively studied.  More
precisely, models over the field of $p$-adic numbers have been
considered which is due to the assumption that $p$-adic numbers
provide a more exact and more adequate description of micro-world
phenomena. Numerous applications of these numbers to theoretical
physics have been proposed in papers \cite{ADFV}, \cite{FW},
\cite{MP}, \cite{V1},\cite{V2} to quantum mechanics \cite{Kh1}, to
$p$-adic - valued physical observable \cite{Kh1} and many others
\cite{Kh2},\cite{VVZ}.

The study of $p$-adic dynamical systems arises in Diophantine
geometry in the constructions of canonical heights, used for
counting rational points on algebraic vertices over a number
field, as in \cite{CS}. In \cite{Kh4},\cite{TVW} $p$-adic field
have arisen in physics in the theory of superstrings, promoting
questions about their dynamics. Also some applications of $p$-adic
dynamical systems to some biological, physical systems has been
proposed in
\cite{ABKO},\cite{AKK},\cite{DGKS},\cite{Kh4},\cite{Kh5}. Other
studies of non-Archimedean dynamics in the neighborhood of a
periodic and of the counting of periodic points over global fields
using local fields appear in \cite{HY},\cite{L},\cite{P}. It is
known that the analytic functions play important role in complex
analysis. In the $p$-adic analysis the rational functions play a
similar role to the analytic functions in complex analysis
\cite{R}. Therefore, naturally one arises a question to study the
dynamics of these functions in the $p$-adic analysis. On the hand,
such $p$-adic dynamical systems appear while studying $p$-adic
Gibbs measures \cite{gmr}. In \cite{B1},\cite{B2} dynamics on the
Fatou set of a rational function defined over some finite
extension of $\Q_p$ have been studied, besides, an analogue of
Sullivan's no wandering domains theorem for $p$-adic rational
functions which have no wild recurrent Julia critical points were
proved. In \cite{AKTS}  the behavior  of a $p$-adic dynamical
system $f(x)=x^n$ in the fields of $p$-adic numbers $\Q_p$ and
complex $p$-adic numbers $\C_p$ was investigated. Some ergodic
properties that dynamical system has been considered in
\cite{GKL}.

The base of $p$-adic analysis, $p$-adic mathematical physics
are explained in \cite{G},\cite{Ko},\cite{VVZ}.

In the paper we will investigate the behavior of trajectory of
rational $p$-adic dynamical systems in $\C_p$. We will study Sigel
disks and attractors of such dynamical systems. In the final
section we show the considered dynamical system is not ergodic.

\section{Preliminaries}

\subsection{$p$-adic numbers}

Let $\Q$ be the field of rational numbers. The greatest common
divisor of the positive integers $n$ and $m$ is denotes by
$(n,m)$. Every rational number $x\neq 0$ can be represented in the
form $x=p^r\dsf{n}{m}$, where $r,n\in\mathbb{Z}$, $m$ is a
positive integer, $(p,n)=1$, $(p,m)=1$ and $p$ is a fixed prime
number. The $p$-adic norm of $x$ is given by
$$
|x|_p=\left\{
\ba{ll}
p^{-r} & \ \textrm{ for $x\neq 0$}\\
0 &\ \textrm{ for $x=0$}.\\
\ea
\right.
$$
It has the following properties:

1) $|x|_p\geq 0$ and $|x|_p=0$ if and only if $x=0$,

2) $|xy|_p=|x|_p|y|_p$,

3) the strong triangle inequality
$$
|x+y|_p\leq\max\{|x|_p,|y|_p\},
$$

3.1) if $|x|_p\neq |y|_p$ then $|x-y|_p=\max\{|x|_p,|y|_p\}$,

3.2) if $|x|_p=|y|_p$ then $|x-y|_p\leq |x|_p$,

this is a non-Archimedean one.

The completion of $\Q$ with  respect to $p$-adic norm defines
the $p$-adic field which is denoted by $\Q_p$.

The well-known Ostrovsky's theorem asserts that norms $|x|=|x|_{\infty}$ and
$|x|_p$, $p=2,3,5...$ exhaust all nonequivalent norms on $\Q$ (see \cite{K}).
Any $p$-adic number $x\neq 0$ can be uniquely represented in the
canonical series:
$$
x=p^{\g(x)}(x_0+x_1p+x_2p^2+...) ,
\eqno(2.1)
$$
where $\g=\g(x)\in\Z$ and $x_j$ are integers, $0\leq x_j\leq p-1$, $x_0>0$,
$j=0,1,2,...$ (see more detail \cite{G},\cite{Ko}).
Observe that in this case $|x|_p=p^{-\g(x)}$.

The algebraic completion of $\Q_p$ is denoted by $\C_p$ and it is called
{\it complex $p$-adic numbers}.  For any $a\in\C_p$ and $r>0$ denote
$$
U_r(a)=\{x\in\C_p : |x-a|_p\leq r\},\ \ V_r(a)=\{x\in\C_p : |x-a|_p< r\},
$$
$$
S_r(a)=\{x\in\C_p : |x-a|_p= r\}.
$$

A function $f:U_r(a)\to\C_p$ is said to be {\it analytic} if it can be represented by
$$
f(x)=\sum_{n=0}^{\infty}f_n(x-a)^n, \ \ \ f_n\in \C_p,
$$ which converges uniformly on the ball $U_r(a)$.

\subsection{Dynamical systems in $\C_p$}

In this section we recall some known facts concerning dynamical systems $(f,U)$ in $\C_p$,
where $f: x\in U\to f(x)\in U$ is an analytic function and $U=U_r(a)$ or $\C_p$.

Now let $f:U\to U$ be an analytic function. Denote $x_n=f^n(x_0)$, where $x_0\in U$ and
$f^n(x)=\underbrace{f\circ\dots\circ f(x)}_n$.

Recall some  the standard terminology of the theory of dynamical
systems (see for example \cite{PJS}). If $f(x_0)=x_0$ then $x_0$
is called a {\it fixed point}. A fixed point $x_0$ is called an
{\it attractor} if there exists a neighborhood $V(x_0)$ of $x_0$
such that for all points $y\in V(x_0)$ it holds
$\lim\limits_{n\to\infty}y_n=x_0$. If $x_0$ is an attractor then
its {\it basin of attraction} is
$$
A(x_0)=\{y\in \C_p :\ y_n\to x_0, \ n\to\infty\}.
$$
A fixed point $x_0$ is called {\it repeller} if there  exists a
neighborhood $V(x_0)$ of $x_0$ such that $|f(x)-x_0|_p>|x-x_0|_p$
for $x\in V(x_0)$, $x\neq x_0$. Let $x_0$ be a fixed point of a
function $f(x)$. The ball $V_r(x_0)$ (contained in $U$) is said to
be a {\it Siegel disk} if each sphere $S_{\r}(x_0)$, $\r<r$ is an
invariant sphere of $f(x)$, i.e. if $x\in S_{\r}(x_0)$ then all
iterated points $x_n\in S_{\r}(x_0)$ for all $n=1,2\dots$.  The
union of all Siegel desks with the center at $x_0$ is said to {\it
a maximum Siegel disk} and is denoted by $SI(x_0)$.

{\bf Remark.}\cite{AKTS} In complex geometry, the center of a disk
is uniquely determined by the disk, and different fixed points
cannot have the same Siegel disks. In non-Archimedean geometry, a
center of a disk is nothing but a point which belongs to the disk.
Therefore, in principle, different fixed points may have the same
Siegel desk.

Let $x_0$ be a fixed point of an analytic function  $f(x)$. Put
$$
\l=\frac{d}{dx}f(x_0).
$$

The point $x_0$ is called {\it attractive} if $0\leq |\l|_p<1$, {\it indifferent} if
$|\l|_p=1$, and {\it repelling} if $|\l|_p>1$.

{\bf Theorem 2.1.}\cite{AKTS} {\it Let $x_0$ be a fixed point of an analytic function $f:U\to U$.
The following assertions hold

1. if $x_0$ is an attractive point of $f$, then it is an attractor
of the dynamical system $(f,U)$. If $r>0$ satisfies the inequality
$$
q=\max_{1\leq n<\infty}\bigg|\frac{1}{n!}\frac{d^nf}{dx^n}(x_0)\bigg|_pr^{n-1}<1
\eqno(2.2)
$$
and $U_r(x_0)\subset U$ then $U_r(x_0)\subset A(x_0)$;

2. if $x_0$ is an indifferent point of $f$ then it is the center of a Siegel disk. If $r$
satisfies the inequality
$$
s=\max_{2\leq n<\infty}\bigg|\frac{1}{n!}\frac{d^nf}{dx^n}(x_0)\bigg|_pr^{n-1}<|f'(x_0)|_p
\eqno(2.3)
$$
and $U_r(x_0)\subset U$ then $U_r(x_0)\subset SI(x_0)$;

3. if $x_0$ is a repelling point of $f$ then $x_0$ is a repeller of the dynamical
system $(f,U)$.}

\section{Rational $p$-adic dynamical systems}

In this section we consider dynamical system associated with the function $f:\C_p\to\C_p$ defined by
$$
f(x)=\frac{x+a}{bx+c}, \ \ b\neq 0,\ c\neq ab, \ a,b,c\in \C_p
\eqno(3.1)
$$
where  $x\neq \hat x=-\dsf{c}{b}$.

It is not difficult to check that fixed points of the function (3.1) are
$$
x_{1,2}=\frac{1-c\pm\sqrt{(c-1)^2+4ab}}{2b}.
\eqno(3.2)
$$

The following theorem is important in our investigation.

{\bf Theorem 3.1.} {\it Let $x_1$ and $x_2$ be the fixed points of (3.1) (see (3.2)). Then

1. the point $x_1$ (resp. $x_2$) is attractive if and only if the point $x_2$
(resp. $x_1$) is  repelling.

2. the point $x_1$ is  indifferent if and only if the point $x_2$ is one.}

The proof immediately follows from the easily checking  equality
$$
f'(x_1)\cdot f'(x_2)=1.
$$

\subsection{Case: $|f'(x_1)|_p=1$}

Let $|f'(x_1)|_p=1$, then according to Theorem 3.1 we have $|f'(x_2)|_p=1$. Observe that
the considered case is equivalent to condition
$$
|f'(x_i)|_p=\bigg|\frac{c-ab}{(bx_1+c)^2}\bigg|_p=1, \ \ i=1,2.
\eqno(3.3)
$$

By Theorem 2.1 every fixed point $x_i$ is a center of Siegel disk. We now verify the
condition (2.3). First of all we compute
\bea
\bigg|\frac{1}{n!}\frac{d^nf}{dx^n}(x_i)\bigg|_p=
|c-ab|_p\bigg|\frac{nb^{n-1}}{(bx_i+c)^{n+1}}\bigg|_p\leq \nonumber \\
\leq \bigg|\bigg(\frac{b}{bx_i+c}\bigg)^{n-1}\bigg|_p= \nonumber \\
=\bigg|\frac{b}{\sqrt{c-ab}}\bigg|^{n-1}_p, \nonumber
\eea
here we have used the equality (3.3). Then the condition (2.3) is
satisfied if the following inequality holds
$$
\max_{2\leq n<\infty}\bigg|\frac{b}{\sqrt{c-ab}}\bigg|^{n-1}_pr^{n-1}<1.
\eqno(3.4)
$$
Let
$$
\bigg|\frac{b}{\sqrt{c-ab}}\bigg|_p<1.
\eqno(3.5)
$$
If $r\leq 1$ then the condition (3.4) is satisfied, and hence $U_1(x_i)\subset SI(x_i)$.

{\bf Theorem 3.2.} {\it Let the conditions (3.3) and (3.5) be satisfied. Then
$$
SI(x_i)=V_{1+\i_c}(x_i), \ \ i=1,2,
$$
where $\i_c=\bigg|\dsf{\sqrt{c-ab}}{b}\bigg|_p-1$.}

{\it Proof.} It suffices to prove that for any $\i<\i_c$ the
equality
$$
f(S_{1+\i}(x_i))=S_{1+\i}(x_i)
\eqno(3.6)
$$
is valid. Let $y\in S_{1+\i}(x_i)$, i.e. $y=x_i+\g$, where $|\g|_p=1+\i$.
Then from (3.1) we get
$$
|f(y)-x_i|_p=\frac{1+\i}{\bigg|\dsf{c-ab}{(bx_i+c)^2}+\dsf{\g b}{bx_i+c}\bigg|_p}.
\eqno(3.7)
$$
The inequality $\i<\i_c$ implies that $\bigg|\dsf{\g b}{bx_i+c}\bigg|_p<1$. It then
follows from (3.7) and (3.3) that $|f(y)-x_i|_p=1+\i$, which means that (3.6) is valid. Here
we have used 3.1) property of the  norm $|\cdot |_p$.
If $\i>\i_c$ then  $\bigg|\dsf{\g b}{bx_i+c}\bigg|_p>1$, consequently from (3.7) we infer
$|f(y)-x_i|_p=1+\i_c<1+\i$. Hence (3.6) does not hold. One remains
to consider the case $\i=\i_c$. We choose
$\g$ as follows
$$
\g=\tilde\g=(p-1)\frac{c-ab}{b(bx_1+c)}.
$$
Then it easy to check that $\tilde y=x_i+\tilde\g$ belongs $S_{1+\i_c}(x_i)$, but
$$
|f(\tilde y)-x_i|_p=p(1+\i_c)>1+\i_c.
$$
Thus the equality (3.6) is valid only at $\i<\i_c$. This completes the proof.

Now we are  interested the relation between Siegel disks $SI(x_i)$, $i=1,2.$

{\bf Theorem 3.3.} {\it  Let the conditions (3.3) and (3.5) be satisfied.

(i) If $\bigg|\dsf{\sqrt{(c-1)^2+4ab}}{b}\bigg|_p\geq 1+\i_c$,
then $SI(x_1)\cap SI(x_2)=\emptyset$;

(ii) otherwise   $SI(x_1)=SI(x_2)$.}

{\it Proof.} (i) From (3.2) we find
$$
|x_1-x_2|_p=\bigg|\frac{\sqrt{(c-1)^2+4ab}}{b}\bigg|_p\geq 1+\i_c.
$$
This means $x_1\notin V_{1+\i_c}(x_2)$, hence by Theorem 3.2 we have
$SI(x_1)\cap SI(x_2)=\emptyset$.

(ii) In this case we have $|x_1-x_2|_p<1+\i_c$. Let $y\in SI(x_1)$, then by
Theorem 3.2 we can write $|y-x_1|_p<1+\i_c$. Whence
$$
|y-x_2|_p=|(y-x_1)+(x_1-x_2)|_p<1+\i_c.
$$
Consequently, we have $SI(x_1)\subset SI(x_2)$. So
$SI(x_1)=SI(x_2)$ since balls with the same radius ether coincide
or disjoint in Non-Archimedean setting.

\subsection{Case: $|f'(x_1)|_p\neq 1$}

According to Theorem 3.1 without loss of generality we may assume that $|f'(x_1)|_p<1$.In this case we have $|f'(x_2)|_p>1$.

From Theorem 2.1 we obtain the following

{\bf Proposition 3.4.} {\it The fixed point $x_2$ is a repelling point of the dynamical
system.}

Now one remains to investigate the fixed point $x_1$. Observe that the condition
$|f'(x_1)|_p<1$ is equivalent to
$$
\bigg|\frac{c-ab}{(bx_1+c)^2}\bigg|_p<1.
\eqno(3.8)
$$
Suppose the following condition to be satisfied
$$
\bigg|\frac{b}{bx_1+c}\bigg|_p=\bigg|\frac{2b}{1+c+\sqrt{(c-1)^2+4ab}}\bigg|_p\leq 1.
\eqno(3.9)
$$

{\bf Lemma 3.5.} {\it Let the conditions (3.8) and (3.9) be satisfied. Then
the inclusion
$$
V_1(x_1)\subset A(x_1)
$$
is valid.}

{\it Proof.} We check the condition (2.2) of Theorem 2.1.
$$
q=\max_{1\leq n<\infty}\bigg|\frac{n(c-ab)b^{n-1}}{(bx_1+c)^{n+1}}\bigg|_pr^{n-1}<
\max_{1\leq n<\infty}\bigg|\frac{b}{bx_1+c}\bigg|^{n-1}_pr^{n-1}<1.
$$
According to (3.9) this condition is fulfilled if $r<1$. By Theorem 2.1 we infer
the required assertion.

Denote
$$
\d_1=\bigg|\frac{(bx_1+c)^2}{c-ab}\bigg|_p-1, \ \ \d_2=\bigg|\frac{bx_1+c}{b}\bigg|_p-1.
$$

{\bf Lemma 3.6.} {\it Let  the conditions (3.8) and (3.9) be satisfied. Then
$\hat x,x_2\in S_{1+\d_2}(x_1)$, here
$$
\hat x=-\frac{c}{b}, \ \ x_2=\frac{1-c-\sqrt{(c-1)^2+4ab}}{2b}.
$$}

{\it Proof.} Consider \bea
|\hat x-x_1|_p=\bigg|\frac{1+c+\sqrt{(c-1)^2+4ab}}{2b}\bigg|_p=\nonumber \\
=\bigg|\frac{bx_1+c}{b}\bigg|_p=1+\d_2, \nonumber
\eea
hence $\hat x\in S_{1+\d_2}(x_1)$. It easy to check that
\bea
|\hat x-x_2|_p=\bigg|\frac{c-ab}{b(bx_1+c)}\bigg|_p=\nonumber\\
=\frac{1+\d_2}{1+\d_1}<1+\d_2, \nonumber
\eea
since $\d_1>0$. So  we have
$$
|x_1-x_2|_p=|(x_1-\hat x) +(\hat x - x_2)|_p=1+\d_2,
$$
here we have used 3.1) property of the norm. Lemma is proved.

{\bf Theorem 3.7.} {\it Let the conditions (3.8) and
$$
\bigg|\frac{b}{bx_1+c}\bigg|_p<1
\eqno(3.10)
$$
be satisfied. Then
$$
\bigcup_{\d: 0\leq \d\neq 1+\d_2}S_{\d}(x_1)\subset A(x_1).
$$}

{\it Proof.} It suffices to prove that $S_{1+\d}(x_1)\subset
A(x_1)$ at $\d\neq 1+\d_2$. Indeed, from the proof of Lemma 3.5
one easily sees that the condition (3.10) provides
$U_1(x_1)\subset A(x_1)$. Let $x\in S_{1+\d}(x_1)$, i.e.
$x=x_1+\g$, $|\g|_p=1+\d$. From (3.1) and (3.2) we get
$$
|f(x)-x_1|_p=\bigg|\frac{(c-ab)\g}{(bx_1+c)^2+\g b(bx_1+c)}\bigg|_p=
$$
$$
=\frac{|c-ab|(1+\d)}{|b|^2\bigg|\bigg(\dsf{bx_1+c}{b}\bigg)^2+\g\bigg(\dsf{bx_1+c}{b}\bigg)\bigg|_p}.
\eqno(3.11)
$$

{\tt Case 1.} Let us assume that
$$
\bigg|\dsf{bx_1+c}{b}\bigg|_p^2>|\g|_p\bigg|\dsf{bx_1+c}{b}\bigg|_p,
$$ this implies that $\d<\d_2$. Then from (3.11) we infer
$$
|f(x)-x_1|_p=\bigg|\frac{c-ab}{(bx_1+c)^2}\bigg|_p(1+\d).
\eqno(3.12)
$$
If $\d\leq\d_1$ then the right side of (3.11) is not greeter  than $1$. Hence, $f(x)\in U_1(x_1)$.
This yields that $S_{1+\d}(x_1)\subset A(x_1)$. If $\d>\d_1$, then the right side
of (3.12) is greeter than 1, denote it by  $1+\l$, i.e.
$$
1+\l=\bigg|\frac{c-ab}{(bx_1+c)^2}\bigg|_p(1+\d), \ \ \l>0.
$$
From (3.8) and $\d<\d_2$ we obtain $\l<\d_2$, since in this case $f(x)\in S_{1+\l}(x_1)$,
so we can put $f(x)$ instead of  $x$ in (3.12), namely
$$
|f^{2}(x)-x_1|_p=\bigg|\frac{c-ab}{(bx_1+c)^2}\bigg|_p|f(x)-x_1|_p=
\bigg|\frac{c-ab}{(bx_1+c)^2}\bigg|^2_p(1+\d) .
$$
If the right side of the last equality is not greeter than 1, then $f^{2}(x)\in U_1(x_1)$,
and hence $S_{1+\d}(x_1)\subset A(x_1)$. Otherwise repeating the above argument
we can prove the following equality
$$
|f^{n}(x)-x_1|_p=\bigg|\frac{c-ab}{(bx_1+c)^2}\bigg|^n_p(1+\d), \ \ n\geq 1.
\eqno(3.13)
$$
The condition (3.8) implies that there is a positive integer $n_0$ such that
$f^{n}(x)\in U_1(x_1)$ for all $n>n_0$. This yields $S_{1+\d}(x_1)\subset A(x_1)$.

{\tt Case 2.} Now suppose that
$$
\bigg|\dsf{bx_1+c}{b}\bigg|_p^2<|\g|_p\bigg|\dsf{bx_1+c}{b}\bigg|_p,
$$
this implies that $\d>\d_2$. It then follows from (3.11) that
$$
|f(x)-x_1|_p=\bigg|\dsf{c-ab}{b(bx_1+c)}\bigg|_p.
\eqno(3.14)
$$
Observe that
$$
\bigg|\dsf{c-ab}{b(bx_1+c)}\bigg|_p=\frac{1+\d_2}{1+\d_1}.
\eqno(3.15)
$$
If $\d_2\leq \d_1$ then the equalities (3.14) and (3.15) provide that $f(x)\in U_1(x_1),$
and in this case we obtain the assertion of theorem. If $\d_2>\d_1$, then from
(3.15) we infer that the right side of (3.14) is greeter than 1, and which is
denoted by $1+\m$, $\m>0$. So $f(x)\in S_{1+\m}(x_1)$. Show that $\m<\d_2$. Indeed,
$$
1+\m=\frac{1+\d_2}{1+\d_1}<1+\d_2
$$
this implies $\m<\d_2$. Thus, we have reduced our consideration to the case 1.
This completes the proof.

{\bf Lemma 3.8.} {\it Let $|x-x_2|_p>\dsf{1+\d_2}{1+\d_1}$, then
$f(x)\in S_{1+\d_2}(x_2)$.}

{\it Proof.} Denote  $\g=x-x_2$ , then the condition of lemma
means that
$$
\bigg|\frac{\g b}{bx_2+c}\bigg|_p>1.
$$
we then have \bea
|f(x)-x_2|_p=\frac{|c-ab|_p|x-x_2|_p}{|(bx_2+c)^2+\g b(bx_2+c)|_p}=\nonumber \\
=\bigg|\dsf{c-ab}{(bx_2+c)^2}\bigg|_p\dsf{|x-x_2|_p}{\bigg|1+\dsf{\g b}{bx_2+c}\bigg|_p}=\nonumber\\
=\bigg|\dsf{c-ab}{(bx_2+c)^2}\bigg|_p\dsf{|x-x_2|_p}{\bigg|\dsf{\g b}{bx_2+c}\bigg|_p}=\nonumber\\
=\bigg|\dsf{c-ab}{b(bx_2+c)}\bigg|_p=\nonumber\\
=\bigg|\dsf{bx_1+c}{b}\bigg|_p=1+\d_2\nonumber
\eea
Lemma is proved.

From this lemma we obtain the following

{\bf Corollary 3.9.} {\it Let the condition of Lemma 3.8 be satisfied, then
$f(S_{1+\d_2}(x_2))=S_{1+\d_2}(x_2)$.}

{\bf Corollary 3.10.} {\it Let $|x-x_2|_p\leq \dsf{1+\d_2}{1+\d_1}$, then
$$
|f(x)-x_2|_p\geq |f'(x_2)|_p|x-x_2|_p
$$
}
The proof is similar to the proof of Lemma 3.8.

{\bf Theorem 3.11.} {\it Let the conditions of Theorem 3.7 be satisfied. Then
$$
A(x_1)=\C_p\setminus\{\hat x,x_2\}.
$$
}

{\it Proof.} Since $\hat x$ does not belong the domain of $f$ and
$x_2$ is a fixed point of one, therefore $\hat x,x_2\notin
A(x_1)$.  According to Theorem 3.7 it suffices to prove that
$S_{1+\d_2}(x_1)\setminus\{\hat x,x_2\}\subset A(x_1)$. Keeping in
mind (3.15), $x_1-\hat x=\dsf{bx_1+c}{b}$ and $\g=x-x_1$ the
equality (3.11) yields
$$
|f(x)-x_1|_p=\frac{1+\d_2}{1+\d_1}\cdot\dsf{|x-x_1|_p}{|x-\hat x|_p}.
\eqno(3.16)
$$
From $|x-x_1|_p=1+\d_2$ we get
$$
|f(x)-x_1|_p=\frac{(1+\d_2)^2}{(1+\d_1)|x-\hat x|_p}.
\eqno(3.17)
$$
If the right side of (3.17) non equal to $1+\d_2$ then according to Theorem 3.7
we infer that $f(x)\in A(x_1)$, hence $x\in A(x_1)$. One remains to consider
a case when the right side of (3.17) is equal to $1+\d_2$, i.e
$|f(x)-x_1|_p=1+\d_2$. In this case we find
$$
|x-\hat x|_p=\dsf{1+\d_2}{1+\d_1}.
\eqno(3.18)
$$
From this we get
$$
|x-x_2|_p=|(x-\hat x)+(\hat x-x_2|_p\leq \dsf{1+\d_2}{1+\d_1},
\eqno(3.19)
$$
here it has been used the equality $|\hat x-x_2|_p=\dsf{1+\d_2}{1+\d_1}$.
(see the proof of Lemma 3.6). According to Corollary 3.10 there  exists
a positive integer $n_0$ such that
$$
|f^{n_0}(x)-x_2|_p>\dsf{1+\d_2}{1+\d_1},
$$
whence
$$
|f^{n_0}(x)-\hat x|_p=|(f^{n_0}(x)-x_2)+(x_2-\hat x)|_p>\dsf{1+\d_2}{1+\d_1},
\eqno(3.20)
$$
here we have used 3.1) property of the norm. From (3.16) we obtain
$$
|f^{2}(x)-x_1|_p=\frac{1+\d_2}{1+\d_1}\cdot\dsf{|f(x)-x_1|_p}{|f(x)-\hat x|_p}=
\frac{(1+\d_2)^2}{1+\d_1}\cdot\dsf{1}{|f(x)-\hat x|_p}.
\eqno(3.21)
$$
Now we estimate $|f(x)-\hat x|_p$:
$$
|f(x)-\hat x|_p=|(f(x)-x_1)+(x_1-\hat x)|_p\leq 1+\d_2.
$$
Hence the equality (3.21) implies
$$
|f^{2}(x)-x_1|_p\geq \frac{1+\d_2}{1+\d_1}.
\eqno(3.22)
$$
If the left side of (3.22) is not equal to $1+\d_2$, then $f^{2}(x)\in A(x_1)$,
so $x\in A(x_1)$. If $|f^{2}(x)-x_1|_p=1+\d_2$, then repeating this argument
till a number $k$ such that
$$
|f^{k}(x)-x_1|_p\neq 1+\d_2, \eqno(3.23)
$$
we conclude that $x\in A(x_1)$. Now we show that such  number $k$
does exist. Let us assume that for all $m\leq n_0$ the following
equality
$$
|f^{m}(x)-x_1|_p=1+\d_2, \eqno(3.24)
$$
is valid. Otherwise nothing to prove. Put $k=n_0+1$. According to
(3.16),(3.20) and (3.24) we have
$$
|f^{n_0+1}(x)-x_1|_p=\frac{1+\d_2}{1+\d_1}\cdot
\dsf{|f^{n_0}(x)-x_1|_p}{|f^{n_0}(x)-\hat x|_p}<1+\d_2,
$$
i.e. (3.23) is valid. This completes the proof.

\section{Dynamical system $f(x)=\dsf{x}{bx+c}$ in $\Q_p$ is not
ergodic}

In this section we consider a dynamical system
$$
f(x)=\dsf{x}{bx+c}, \ \ c\neq 0,\ b,c\in\Q_p \eqno(4.1)
$$ in $\Q_p$. It is easy to that $x=0$ is a fixed point for (4.1).
A question about ergodicity of the considered system arises when
the fixed point $x=0$ is indifferent. This lead us to the
condition $|c|_p=1$. From condition (3.5) we find that
$|b|_p<|c|_p$.  Then it is not difficult to check that
$f(S_1(0))=S_1(0)$. From now we consider the dynamical system
(4.1) on the sphere $S_1(0)$.

{\bf Lemma 4.1.}{\it  For every ball $U_{p^{-l}}(a)\subset S_1(0)$
then the following equality holds
$$
f(U_{p^{-l}}(a))=U_{p^{-l}}(f(a))
$$
}

{\it Proof.} From inclusion $U_{p^{-l}}(a)\subset S_1(0)$ we have
$|a|_p=1$. Let $|x-a|_p\leq p^{-l}$, then
$$
|f(x)-f(a)|_p=\frac{|c|_p|x-a|_p}{|bx+c|_p|ba+c|_p}=|x-a|\leq
p^{-l},
$$
here we have used the equality $|bx+c|_p=1$, which follows from
$|b|_p<1$. Lemma is proved.

Consider a measurable space $(S_1(0),{\cal B})$, here ${\cal B}$
is the algebra of generated by clopen subsets of $S_1(0)$. Every
element of ${\cal B}$ is a union of some balls $U_{p^{-l}}(a)$. A
measure $\mu:{\cal B}\to \R$ is said to be {\it Haar measure} if
it is defined by
$$
\mu(U_{p^{-l}}(a))=\frac{1}{q^l},
$$
here $q$ is a prime number.

From lemma 4.1 we conclude that $f$ preserves the measure $\mu$,
i.e.
$$
\mu(f(U_{p^{-l}}(a)))=\mu(U_{p^{-l}}(a))\eqno(4.2).$$

Recall a dynamical system $(X,T,\lambda)$, where $T:X\to X$ is a
measure preserving transformation, is called {\it ergodic} if for
every invariant set $A$, i.e. $T(A)=A$ the equalities
$\lambda(A)=0$ or $\lambda(A)=1$ are valid. (see, \cite{W})

{\bf Proposition 4.2.} {\it If there is some number $N\in\N$ such
that $|f^N(a)-a|_p<1$ for some $a\in S_1(0)$ then the dynamical
system (4.1) is not ergodic on $S_1(0)$ with respect to the Haar
measure.}

{\it Proof.} Denote $N=\min\{n\in\N : |f^n(a)-a|_p<1\ \ \textrm{
for some}\ \ a\in S_1(0)\}$. Because of the discreteness of the
$p$-adic metric we can assume that $|f^N(a)-a|_p\leq p^{-l}$ for
some positive integer $l\in\N$. Put
$$
A=\bigcup_{k=0}^{N-1}U_{p^{-l}}(f^k(a)).
$$
Then from Lemma 4.1 we find that $f(A)=A$. It is clear that
$\mu(A)\neq 0$ and $\mu(S_1(0)\setminus A)\neq 0$, hence $f$ is
not ergodic. This completes the proof.\\

{\bf Corollary 4.3.} {\it If $p=2$ then the dynamical system (4.1)
is not ergodic on $S_1(0)$ with respect to the Haar measure.}

{\it Proof.} Form the condition $|c|_2=1$ using the property 3.2)
of the norm we get $|1-c|_2\leq\dsf{1}{2}$. Then we have
$$
|f(a)-a|_2=|1-c-ba|_2\leq\frac{1}{2},
$$
since $|ba|_2\leq\dsf{1}{2}$. Hence the set $A=U_{2^{-1}}(a)$ is
invariant with respect to $f$. On the other hand $\mu(A)=1/2$,
that means $f$ is not ergodic. The corollary is proved.

{\bf Lemma 4.4.}{\it  For every $N\in\N$ the following equalities
hold
$$
f^2(x)=\frac{x}{bx(1+c)+c^{2}},  \ \
f^3(x)=\frac{x}{bx(1+c^2+c^3)+c^{4}}, \ \
f^N(x)=\frac{x}{bx(1+S_N)+c^{2^{N-1}}},
$$
where $S_N=\sum\limits_{m=0}^{N-3}c^{2^m(2^{N-2}-1)}$.}

The proof immediately follows from induction method.\\
From this Lemma we can prove the following

{\bf Corollary 4.5.}{\it  For every integer $N\in\N$ the following
equalities hold
$$f^{-1}(x)=\frac{x}{-bx/c+1/c},  \ \
f^{-2}(x)=\frac{x}{-bx/c(1+1/c)+1/c^{2}}, $$ $$ \ \
f^{-3}(x)=\frac{x}{-bx/c(1+1/c^2+1/c^3)+1/c^{4}},
f^{-N}(x)=\frac{x}{-bx/c(1+Z_N)+1/c^{2^{N-1}}},
$$
where $Z_N=\sum\limits_{m=0}^{N-3}1/c^{2^m(2^{N-2}-1)}$.}\\

Using Lemma 4.4 consider the difference
$$
|f^N(x)-x|_p=|1-bx(1+S_N)-c^{2^{N-1}}|_p
$$
form this we conclude that the condition $|f^N(x)-x|_p=1$ for all
$N\in\N$ and $x\in S_1(0)$ is equivalent to the equality
$$
|1-c^{2^{N-1}}|_p=1, \ \ \forall N\in\N .\eqno(4.3)
$$

Using Corollary 4.5 and analogous argument as above we can obtain
that the condition $|f^{-N}(x)-x|_p=1$ for all $N\in\N$ and $x\in
S_1(0)$ is also equivalent to the equality (4.3).

From (2.1) we have that every element $c\in S_1(0)$ is represented
in the form
$$
c=a_0+a_1p+a_2p^2+...,
$$
where $a_0\neq 0,\ a_k\in\{0,1,...,p-1\},\ k\in\N$. Then it is
easy to see the condition (4.3) is  equivalent to the following
one
$$
a_0^{2^{N-1}}\equiv\!\!\!\!\!\! / \ \ 1(\textrm{mod $p$}), \ \
\forall
N\in\N
$$
This condition is satisfied for example on $p=7$ with $a_0=2$.

Let us assume that the condition (4.2) is satisfied. Then for all
$n\in\N$, we get
$$U_{p^{-l}}(f^{-n}(a))\cap U_{p^{-l}}(a)=\emptyset$$
for all $a\in S_1(0)$,\ $l\in\N$. Then according to Theorem
1.5\cite{W} we conclude that the set the dynamical system (4.1) is
not ergodic.  So we have proved the following

{\bf Theorem 4.6.} {\it The dynamical system (4.1) is not ergodic
on $S_1(0)$ with respect to the Haar measure.}

{\bf Acknowledgement}

This work was done in the scheme of Mathematical Fellowship at the
Abdus Salam International Center for Theoretical Physics (ICTP),
and the authors thank ICTP for providing financial support and all
facilities. The first author also thanks IMU/CDE- program for
financial support. The authors also acknowledge with gratitude to
Professor A.Yu.Khrennikov for the helpful comments.

\end{document}